# SOME STRONG LIMIT THEOREMS FOR THE LARGEST ENTRIES OF SAMPLE CORRELATION MATRICES

By Deli Li[1] and Andrew Rosalsky

*Lakehead University and University of Florida*

Let $\{X_{k,i}; i \ge 1, k \ge 1\}$ be an array of i.i.d. random variables and let $\{p_n; n \ge 1\}$ be a sequence of positive integers such that $n/p_n$ is bounded away from 0 and $\infty$. For $W_n = \max_{1 \le i < j \le p_n} |\sum_{k=1}^n X_{k,i}X_{k,j}|$ and $L_n = \max_{1 \le i < j \le p_n} |\hat{\rho}_{i,j}^{(n)}|$ where $\hat{\rho}_{i,j}^{(n)}$ denotes the Pearson correlation coefficient between $(X_{1,i}, \dots, X_{n,i})'$ and $(X_{1,j}, \dots, X_{n,j})'$, the limit laws (i) $\lim_{n\to\infty} \frac{W_n}{n^\alpha} = 0$ a.s. $(\alpha > 1/2)$, (ii) $\lim_{n\to\infty} n^{1-\alpha} L_n = 0$ a.s. $(1/2 < \alpha \le 1)$, (iii) $\lim_{n\to\infty} \frac{W_n}{\sqrt{n \log n}} = 2$ a.s. and (iv) $\lim_{n\to\infty} (\frac{n}{\log n})^{1/2} L_n = 2$ a.s. are shown to hold under optimal sets of conditions. These results follow from some general theorems proved for arrays of i.i.d. two-dimensional random vectors. The converses of the limit laws (i) and (iii) are also established. The current work was inspired by Jiang's study of the asymptotic behavior of the largest entries of sample correlation matrices.

**1. Introduction.** At the origin of the current investigation is the statistical hypothesis testing problem studied by Jiang [7] using the asymptotic distribution of the largest entry of a sample correlation matrix. Jiang's [7] work will now be discussed. Consider a $p$-variate population $(p \ge 2)$ represented by a random vector $\mathbf{X} = (X_1, \dots, X_p)$ with unknown mean $\boldsymbol{\mu} = (\mu_1, \dots, \mu_p)$, unknown covariance matrix $\boldsymbol{\Sigma}$ and unknown correlation coefficient matrix $\mathbf{R}$. Let $\mathcal{M}_{n,p} = (X_{k,i})_{1 \le k \le n, 1 \le i \le p}$ be an $n \times p$ matrix whose rows are an observed random sample of size $n$ from the $\mathbf{X}$ population; that is, the rows of $\mathcal{M}_{n,p}$ are independent copies of $\mathbf{X}$. Jiang [7] assumed that both $n$ and $p$ are large; more precisely, Jiang [7] assumed that, for some $0 < \gamma < \infty$, $p = p_n \sim \gamma n$ as $n \to \infty$. In contradistinction to classical multivariate data

Received March 2005; revised August 2005.
[1]Supported by a grant from the Natural Sciences and Engineering Research Council of Canada.

analysis wherein the dimension $p$ is fixed (see, e.g., Anderson [2]), in contemporary multivariate data analysis the dimension $p$ can be very large and can vary with $n$ and be comparable with $n$ (see, e.g., Donoho [5] and Johnstone [8] who fittingly illustrate this point with many examples comprising a diversity of applications).

When both $n$ and $p$ are large, Jiang [7] considered the statistical test with null hypothesis $H_0 : \mathbf{R} = \mathbf{I}$, where $\mathbf{I}$ is the $p \times p$ identity matrix. In general, this null hypothesis asserts that the components of $\mathbf{X} = (X_1, \ldots, X_p)$ are uncorrelated whereas when $\mathbf{X}$ has a $p$-variate normal distribution, this null hypothesis asserts that these components are independent.

Jiang's [7] test statistic is extremely intuitive and will now be described. Let $n \geq 2$. Set $\bar{X}_i^{(n)} = \sum_{k=1}^n X_{k,i}/n, 1 \leq i \leq p$. Let $\mathbf{X}_i^{(n)}$ denote the $i$th column of $\mathcal{M}_{n,p}, 1 \leq i \leq p$, let $\mathbf{e} = (1, \ldots, 1)' \in \mathcal{R}^n$ and let $\|\cdot\|$ be the Euclidean norm in $\mathcal{R}^n$. Jiang's [7] test is based on the test statistic

$$L_n = \max_{1 \leq i < j \leq p} |\hat{\rho}_{i,j}^{(n)}|, \tag{1.1}$$

where

$$\begin{aligned} \hat{\rho}_{i,j}^{(n)} &= \frac{\sum_{k=1}^n (X_{k,i} - \bar{X}_i^{(n)})(X_{k,j} - \bar{X}_j^{(n)})}{(\sum_{k=1}^n (X_{k,i} - \bar{X}_i^{(n)})^2)^{1/2} (\sum_{k=1}^n (X_{k,j} - \bar{X}_j^{(n)})^2)^{1/2}} \\ &= \frac{(\mathbf{X}_i^{(n)} - \bar{X}_i^{(n)}\mathbf{e})'(\mathbf{X}_j^{(n)} - \bar{X}_j^{(n)}\mathbf{e})}{\|\mathbf{X}_i^{(n)} - \bar{X}_i^{(n)}\mathbf{e}\| \cdot \|\mathbf{X}_j^{(n)} - \bar{X}_j^{(n)}\mathbf{e}\|} \end{aligned} \tag{1.2}$$

is the Pearson correlation coefficient between the $i$th and $j$th columns of $\mathcal{M}_{n,p}$.

Jiang [7] proved the following two limit theorems concerning the test statistic $L_n$ when $p = p_n$ and $\mathcal{M} = \{X_{k,i}; i \geq 1, k \geq 1\}$ is an array of independent and identically distributed (i.i.d.) random variables. Theorem 1.1 provides a law of the logarithm (LL) for $L_n$ and Theorem 1.2 establishes the asymptotic distribution of $L_n$. The limiting distribution in Theorem 1.2 is a type-I extreme value distribution.

Theorem 1.1 ([7]). *Suppose that $E|X_{1,1}|^r < \infty$ for all $0 < r < 30$. If $\lim_{n\to\infty} \frac{n}{p_n} = \gamma \in (0, \infty)$, then*

$$\lim_{n\to\infty} \left(\frac{n}{\log n}\right)^{1/2} L_n = 2 \qquad \textit{almost surely (a.s.).}$$

Theorem 1.2 ([7]). *Suppose that $E|X_{1,1}|^r < \infty$ for some $r > 30$. If $\lim_{n\to\infty} \frac{n}{p_n} = \gamma \in (0, \infty)$, then*

$$\lim_{n\to\infty} P(nL_n^2 - 4\log n + \log\log n \leq t) = \exp\left\{-\frac{1}{\gamma^2\sqrt{8\pi}} e^{-t/2}\right\}, \quad -\infty < t < \infty.$$

Let $\mathbf{R}_n = (\hat{\rho}_{i,j}^{(n)})_{1 \le i,j \le p_n}$ be the $p_n \times p_n$ sample correlation matrix obtained from $\mathcal{M}_{n,p_n} = (\mathbf{X}_1^{(n)}, \ldots, \mathbf{X}_{p_n}^{(n)})$. As was discussed by Jiang [7], by shifting and scaling each column $\mathbf{X}_i^{(n)}$ of $\mathcal{M}_{n,p_n}$, the new data matrix and $\mathcal{M}_{n,p_n}$ have the same sample correlation matrix $\mathbf{R}_n$. Thus if the population is $p_n$-variate normal, under the null hypothesis that the $p_n$ components of $\mathbf{X}$ are independent, the distribution of $\mathbf{R}_n$ is the same as that generated by a data matrix composed of i.i.d. $N(0,1)$ random variables. Jiang [7] thus obtained the following corollary of Theorems 1.1 and 1.2.

Corollary 1.1 ([7]). *Let $\mathcal{M} = \{X_{k,i}; i \ge 1, k \ge 1\}$ be an array of independent random variables where $X_{k,i} \sim N(\mu_i, \sigma_i^2), i \ge 1, k \ge 1$ where $\sigma_i^2 > 0$, $i \ge 1$. Let the sample correlation matrix $\mathbf{R}_n$ be obtained from $\mathcal{M}_{n,p_n} = \{X_{k,i}; 1 \le i \le p_n, 1 \le k \le n\}, n \ge 1$ where $\{p_n; n \ge 1\}$ is a sequence of positive integers satisfying $\lim_{n\to\infty} \frac{n}{p_n} = \gamma \in (0, \infty)$. Then the conclusions of Theorems 1.1 and 1.2 prevail.*

In the current work, the main results are Kolmogorov–Marcinkiewicz–Zygmund-type strong laws of large numbers (SLLNs) (Theorems 2.1 and 2.2) as well as LLs (Theorems 2.3 and 2.4) for both $\{W_n; n \ge 1\}$ and $\{L_n; n \ge 1\}$ where

$$W_n = \max_{1 \le i < j \le p_n} \left| \sum_{k=1}^{n} X_{k,i} X_{k,j} \right|, \qquad n \ge 1. \tag{1.3}$$

Note that $\sum_{k=1}^{n} X_{k,i} X_{k,j}$ is the $(i,j)$th entry of $\mathcal{M}_{n,p_n}' \mathcal{M}_{n,p_n}$, $1 \le i, j \le n, n \ge 1$. In Theorems 2.1 and 2.3 the conditions are also shown to be necessary. (As in Theorems 1.1 and 1.2, the array $\mathcal{M} = \{X_{k,i}; i \ge 1, k \ge 1\}$ is composed of i.i.d. random variables.) We prove in Theorem 2.4 that Theorem 1.1 holds under substantially weaker moment conditions and the condition $\lim_{n\to\infty} \frac{n}{p_n} = \gamma \in (0, \infty)$ is weakened as well. More specifically, the hypotheses of Theorem 2.4 will be satisfied if $E|X_{1,1}|^6 < \infty$ and $n/p_n$ is bounded away from 0 and $\infty$; $\lim_{n\to\infty} \frac{n}{p_n}$ does not need to exist.

The main tools employed by Jiang [7] in proving Theorem 1.1 are (i) a result of Amosova [1] on probabilities of moderate deviations which sharpens a result of Rubin and Sethuraman [12], and (ii) a special case of Theorem 1 of [3] which is, in turn, a special case of the Chen–Stein Poisson approximation method. In the current work, we use quite a few results from classical probability theory and a recent generalization of the Hoffmann–Jørgensen [6] inequalities due to Li and Rosalsky [10].

The plan of the paper is as follows. Theorems 2.1–2.4 will be stated in Section 2 but their proofs will be deferred until Section 4. In Section 3, three very general results (Theorems 3.1–3.3) will be established concerning arrays of i.i.d. two-dimensional random vectors. These results are of interest in their own right but the last two of them will be used in Section 4 to prove Theorems 2.1–2.4.

**2. The main results.** Throughout, let $\mathcal{M} = \{X_{k,i}; i \geq 1, k \geq 1\}$ be an array of i.i.d. random variables, let $\{p_n; n \geq 1\}$ be a sequence of positive integers, and for $n \geq 1$, consider the $n \times p_n$ matrix $\mathcal{M}_{n,p_n}$ as defined in Section 1. Let $\hat{\rho}_{i,j}^{(n)}$ be defined as in (1.2), $1 \leq i, j \leq n, n \geq 1$. Let $\{L_n; n \geq 1\}$ be as in (1.1) with $p = p_n$ and let $\{W_n; n \geq 1\}$ be as in (1.3).

The first theorem is a Kolmogorov–Marcinkiewicz–Zygmund-type SLLN for $\{W_n; n \geq 1\}$.

THEOREM 2.1. *Suppose that $n/p_n$ is bounded away from 0 and $\infty$. Let $\alpha > 1/2$. Then*

$$\lim_{n\to\infty} \frac{W_n}{n^\alpha} = 0 \qquad a.s. \tag{2.1}$$

*if and only if*

$$\sum_{n=1}^{\infty} P\left(\max_{1\leq i<j\leq n} |X_i X_j| \geq n^\alpha\right) < \infty \tag{2.2}$$

*and*

$$EX_1 = 0 \qquad \text{whenever } \alpha \leq 1.$$

*Here and below $X_i = X_{1,i}, i \geq 1$.*

REMARK 2.1. Note that

$$P\left(\max_{1\leq i<j\leq n} |X_i X_j| \geq n^\alpha\right) \leq n^2 P(|X_1 X_2| \geq n^\alpha), \qquad n \geq 1$$

and so (2.2) holds if

$$\sum_{n=1}^{\infty} n^2 P(|X_1 X_2| \geq n^\alpha) < \infty$$

which is equivalent to $E|X_1 X_2|^{3/\alpha} < \infty$. Thus, (2.2) holds if $E|X_1|^{3/\alpha} < \infty$. Also note that

$$\begin{aligned} P\left(\max_{1\leq i<j\leq n} |X_i X_j| \geq n^\alpha\right) &\wedge \frac{n}{2} P(|X_1 X_2| \geq n^\alpha) \\ &\geq P\left(\max_{1\leq i\leq n/2} |X_i X_{[n/2]+i}| \geq n^\alpha\right) \\ &\geq 1 - \exp\left\{-\frac{n}{2} P(|X_1 X_2| \geq n^\alpha)\right\}, \qquad n \geq 2 \end{aligned}$$

and hence since $1-e^{-x}\sim x$ as $x\to 0$,

$$\sum_{n=2}^{\infty} P\Big(\max_{1\le i\le n/2}|X_iX_{[n/2]+i}|\ge n^{\alpha}\Big)<\infty$$

if and only if

$$\sum_{n=1}^{\infty} nP(|X_1X_2|\ge n^{\alpha})<\infty$$

which is equivalent to $E|X_1X_2|^{2/\alpha}<\infty$. Thus $E|X_1|^{2/\alpha}<\infty$ if (2.2) holds.

The second theorem is a Kolmogorov–Marcinkiewicz–Zygmund-type SLLN for $\{L_n; n\ge 1\}$.

THEOREM 2.2. *Suppose that $n/p_n$ is bounded away from 0 and $\infty$. Let $1/2<\alpha\le 1$. If $X_{1,1}$ is nondegenerate and (2.2) holds, then*

$$\lim_{n\to\infty} n^{1-\alpha}L_n = 0 \qquad a.s. \tag{2.3}$$

The third theorem establishes a LL for $\{W_n; n\ge 1\}$.

THEOREM 2.3. *Suppose that $n/p_n$ is bounded away from 0 and $\infty$. Then*

$$\lim_{n\to\infty}\frac{W_n}{\sqrt{n\log n}} = 2 \qquad a.s. \tag{2.4}$$

*if and only if*

$$EX_1=0, EX_1^2=1 \quad \text{and} \quad \sum_{n=1}^{\infty} P\Big(\max_{1\le i<j\le n}|X_iX_j|\ge\sqrt{n\log n}\Big)<\infty. \tag{2.5}$$

REMARK 2.2. By an argument similar to that in Remark 2.1, the condition

$$\sum_{n=1}^{\infty} P\Big(\max_{1\le i<j\le n}|X_iX_j|\ge\sqrt{n\log n}\Big)<\infty \tag{2.6}$$

is weaker than the condition

$$E\Big(\frac{(X_1X_2)^6}{(\log(e+|X_1X_2|))^3}\Big)<\infty \tag{2.7}$$

but stronger than the condition

$$E\Big(\frac{(X_1X_2)^4}{(\log(e+|X_1X_2|))^2}\Big)<\infty. \tag{2.8}$$

Let $h_1(t) = E(X_1^6 I(|X_1| \le t))$ and $h_2(t) = E(X_1^4 I(|X_1| \le t)), t \ge 0$. Then, by using Fubini's theorem, we can see that (2.7) and (2.8) are, respectively, equivalent to

$$E\left(\frac{X_1^6 h_1(|X_1|)}{(\log(e+|X_1|))^3}\right) < \infty$$

and

$$E\left(\frac{X_1^4 h_2(|X_1|)}{(\log(e+|X_1|))^2}\right) < \infty.$$

Remark 2.3. Note that

$$P\left(\max_{1\le i<j\le n} |X_i X_j| \ge \sqrt{n\log n}\right) = P(Z_{n:1} Z_{n:2} \ge \sqrt{n\log n}\,),$$

where $Z_{n:1}$ and $Z_{n:2}$ are, respectively, the largest and the second largest of the random variables $|X_1|, |X_2|, \ldots, |X_n|$. Thus (2.6) is equivalent to

$$\sum_{n=1}^{\infty} P(Z_{n:1} Z_{n:2} \ge \sqrt{n\log n}) < \infty. \tag{2.9}$$

Clearly,

$$\begin{aligned} P(Z_{n:1} Z_{n:2} \ge \sqrt{n\log n}) &\ge P(Z_{n:2}^2 \ge \sqrt{n\log n}\,) \\ &= P(Z_{n:2} \ge (n\log n)^{1/4}), \qquad n \ge 1. \end{aligned}$$

Let $t_n = P(|X_1| \ge (n\log n)^{1/4}), n \ge 1$. Then (2.6) implies that

$$\begin{aligned} &\sum_{n=1}^{\infty} P(Z_{n:2} \ge (n\log n)^{1/4}) \\ &\qquad = \sum_{n=1}^{\infty} (1 - (1-t_n)^n - nt_n(1-t_n)^{n-1}) \\ &\qquad < \infty \end{aligned}$$

and hence (2.6) implies that $nt_n = o(1)$. These two consequences of (2.6) entail

$$\sum_{n=1}^{\infty} n^2 t_n^2 = \sum_{n=1}^{\infty} n^2 (P(|X_1| \ge (n\log n)^{1/4}))^2 < \infty. \tag{2.10}$$

It follows from the Cauchy–Schwarz inequality and (2.10) that

$$\begin{aligned} &\sum_{n=1}^{\infty} \frac{n^{1/2}}{\log(n+1)} P(X_1^4 > n\log n) \\ &\qquad \le \left(\sum_{n=1}^{\infty} \frac{1}{n\log^2(n+1)}\right)^{1/2} \left(\sum_{n=1}^{\infty} n^2 (P(|X_1| \ge (n\log n)^{1/4}))^2\right)^{1/2} < \infty \end{aligned}$$

and hence (2.6) ensures that

$$E|X_1|^{\beta} < \infty \qquad \text{for all } 0 < \beta < 6.$$

The fourth theorem establishes a LL for $\{L_n; n \geq 1\}$.

THEOREM 2.4. *Suppose that $n/p_n$ is bounded away from 0 and $\infty$. If $X_1$ is nondegenerate and (2.6) holds, then*

$$\lim_{n\to\infty} \left(\frac{n}{\log n}\right)^{1/2} L_n = 2 \qquad a.s. \tag{2.11}$$

**3. Three general results.** In this section three very general results will be established. Theorems 3.2 and 3.3 will be used in Section 4 to prove the four theorems in Section 2.

Let $\{(U_{k,i}, V_{k,i}); i \geq 1, k \geq 1\}$ be an array of i.i.d. two-dimensional random vectors. Let $\{p_n; n \geq 1\}$ be a sequence of positive integers and consider the $n \times p_n$ matrices

$$\mathbf{A}_n \equiv (U_{k,i})_{1 \leq k \leq n, 1 \leq i \leq p_n}, \qquad \mathbf{B}_n \equiv (V_{k,i})_{1 \leq k \leq n, 1 \leq i \leq p_n}, \qquad n \geq 1.$$

Then $\mathbf{A}_n'\mathbf{B}_n$ is a $p_n \times p_n$ matrix whose $(i,j)$th entry is $\sum_{k=1}^{n} U_{k,i}V_{k,j}, n \geq 1$. Let

$$T_n = \max_{1 \leq i \neq j \leq p_n} \left| \sum_{k=1}^{n} U_{k,i}V_{k,j} \right|, \qquad n \geq 1.$$

Let $\{Y_n; n \geq 1\}$ be a sequence of i.i.d. random variables where $Y_1$ has the same distribution as $U_{1,1}V_{1,2}$ and set $S_n = \sum_{k=1}^{n} Y_k, n \geq 1$.

Theorem 3.1 may now be presented. It is a very general result wherein the asymptotic fluctuation behavior of $T_n$ is governed by a Baum–Katz–Lai-type complete convergence result. It is not assumed that $n/p_n$ is bounded away from 0 and $\infty$.

THEOREM 3.1. *Let $\{a_n; n \geq 1\}$ be a sequence of positive constants such that $a_n \uparrow \infty$ and*

$$\lim_{c \downarrow 1} \limsup_{n\to\infty} \frac{a_{[cn]}}{a_n} = 1. \tag{3.1}$$

*Suppose that the sequence $\{p_n; n \geq 1\}$ is nondecreasing. If*

$$\frac{S_n}{a_n} \xrightarrow{P} 0 \tag{3.2}$$

*and*

$$\sum_{n=1}^{\infty} \frac{p_n^2}{n} P\left(\frac{|S_n|}{a_n} > \lambda\right) < \infty \qquad \textit{for some } 0 < \lambda < \infty, \tag{3.3}$$

*then*

$$\limsup_{n\to\infty} \frac{T_n}{a_n} \le \lambda \qquad a.s. \tag{3.4}$$

PROOF. Let $\delta > 0$ be arbitrary. By (3.1), we can choose $c > 1$ such that

$$a_{[cn]} \le (1+\delta)a_n \qquad \text{for all large } n. \tag{3.5}$$

Note that for all large $n$

$$\begin{aligned}\max_{c^{n-1}<m\le c^n} T_m &= \max_{c^{n-1}<m\le c^n} \max_{1\le i\ne j\le p_m} \left|\sum_{k=1}^{m} U_{k,i}V_{k,j}\right| \\ &\le \max_{1\le i\ne j\le p_{[c^n]}} \max_{c^{n-1}<m\le c^n} \left|\sum_{k=1}^{m} U_{k,i}V_{k,j}\right| \\ &\equiv H_n \qquad \text{(say)}\end{aligned}$$

and hence for all large $n$

$$\begin{aligned}&P\left(\max_{c^{n-1}<m\le c^n} \frac{T_m}{a_m} > (1+3\delta)^2\lambda\right) \\ &\quad\le P\left(\frac{H_n}{a_{[c^n]}} > (1+3\delta)\lambda\right) \qquad \text{[by (3.5)]} \\ &\quad\le (p_{[c^n]})^2 P\left(\max_{c^{n-1}<m\le c^n} \frac{|S_m|}{a_{[c^n]}} > (1+3\delta)\lambda\right).\end{aligned}$$

Note that (3.2) ensures that

$$\lim_{n\to\infty} \min_{1\le k\le n} P(S_n - S_k > -\delta\lambda a_n) = 1$$

and

$$\lim_{n\to\infty} \min_{1\le k\le n} P(S_n - S_k < \delta\lambda a_n) = 1.$$

It then follows from Theorem 2.3 of [11] that for all large $n$

$$\begin{aligned}&P\left(\max_{c^{n-1}<m\le c^n} \frac{T_m}{a_m} > (1+3\delta)^2\lambda\right) \\ &\quad\le 2(p_{[c^n]})^2 P\left(\frac{|S_{[c^n]}|}{a_{[c^n]}} > (1+2\delta)\lambda\right).\end{aligned} \tag{3.6}$$

Again recalling (3.2), for all large $n$ and $m \in [[c^n], [c^{n+1}] - 1]$, another application of Theorem 2.3 of [11] yields

$$P\left(\frac{|S_{[c^n]}|}{a_{[c^n]}} > (1+2\delta)\lambda\right)$$

$$(3.7)\qquad \begin{aligned} &\le P\left(\max_{[c^n]\le j\le m}\frac{|S_j|}{a_{[c^n]}}>(1+2\delta)\lambda\right)\\ &\le 2P\left(\frac{|S_m|}{a_{[c^n]}}>(1+\delta)\lambda\right)\\ &\le 2P\left(\frac{(1+\delta)|S_m|}{a_{[c^{n+1}]}}>(1+\delta)\lambda\right)\qquad [\text{by (3.5)}]\\ &\le 2P\left(\frac{|S_m|}{a_m}>\lambda\right). \end{aligned}$$

Since $[c^{n+1}]-[c^n]\sim\frac{c-1}{c}([c^{n+1}]-1)$, it follows from (3.6) and (3.7) that for all large $n$

$$\begin{aligned} &P\left(\max_{c^{n-1}<m\le c^n}\frac{T_m}{a_m}>(1+3\delta)^2\lambda\right)\\ &\quad\le 4(p_{[c^n]})^2\frac{\sum_{m=[c^n]}^{[c^{n+1}]-1}P(|S_m|/a_m>\lambda)}{[c^{n+1}]-[c^n]}\\ &\quad\le \frac{8c}{c-1}(p_{[c^n]})^2\frac{\sum_{m=[c^n]}^{[c^{n+1}]-1}P(|S_m|/a_m>\lambda)}{[c^{n+1}]-1}\\ &\quad\le \frac{8c}{c-1}\sum_{m=[c^n]}^{[c^{n+1}]-1}\frac{p_m^2}{m}P\left(\frac{|S_m|}{a_m}>\lambda\right). \end{aligned}$$

Then by (3.3) and the Borel–Cantelli lemma,

$$P\left(\max_{c^{n-1}<m\le c^n}\frac{T_m}{a_m}>(1+3\delta)^2\lambda\text{ i.o. }(n)\right)=0$$

whence

$$\limsup_{n\to\infty}\frac{T_n}{a_n}\le(1+3\delta)^2\lambda\qquad\text{a.s.}$$

Since $\delta>0$ is arbitrary, the conclusion (3.4) is established. □

Consider the sequence of partial sums $\{S_n;n\ge1\}$ defined prior to the statement of Theorem 3.1. Let $\beta>0$ and $\alpha>1/2$, and assume that $EY_1=0$ if $\alpha\le1$. According to the celebrated theorem of Baum and Katz [4], the following are equivalent:

$$\sum_{n=1}^{\infty}n^{2\beta-1}P\left(\frac{|S_n|}{n^\alpha}>\varepsilon\right)<\infty\qquad\text{for all }\varepsilon>0,$$

$$\sum_{n=1}^{\infty} n^{2\beta-1} P\left(\sup_{m\geq n} \frac{|S_m|}{m^\alpha} > \varepsilon\right) < \infty \qquad \text{for all } \varepsilon > 0,$$

$$E|Y_1|^{(2\beta+1)/\alpha} < \infty. \tag{3.8}$$

Note that (3.8) is equivalent to

$$E|U_{1,1}|^{(2\beta+1)/\alpha} < \infty \quad \text{and} \quad E|V_{1,1}|^{(2\beta+1)/\alpha} < \infty \tag{3.9}$$

and that $EY_1 = 0$ is equivalent to

$$(EU_{1,1})(EV_{1,1}) = 0. \tag{3.10}$$

Combining Theorem 3.1 and the Baum–Katz [4] theorem yields the following.

COROLLARY 3.1. *Let $\alpha > 1/2$ and $\beta > 0$. Suppose that (3.9) holds and if $\alpha \leq 1$ that (3.10) holds. Then*

$$\lim_{n\to\infty} \frac{\max_{1\leq i\neq j\leq n^\beta} |\sum_{k=1}^n U_{k,i}V_{k,j}|}{n^\alpha} = 0 \qquad a.s.$$

PROOF. Let $a_n = n^\alpha$ and $p_n = [n^\beta], n \geq 1$. Then (3.1) is immediate and (3.3) holds for all $\lambda > 0$ by the Baum–Katz [4] theorem. It follows from (3.8) that $E|Y_1|^{1/\alpha} < \infty$ whence by the Kolmogorov–Marcinkiewicz–Zygmund SLLN, (3.2) holds. Thus (3.4) holds for all $\lambda > 0$ by Theorem 3.1. Since $\lambda > 0$ is arbitrary, the corollary is proved. □

Again consider the sequence of partial sums $\{S_n; n \geq 1\}$ defined prior to the statement of Theorem 3.1 and let $\beta > 0$. By a theorem of Lai [9],

$$\sum_{n=2}^{\infty} n^{2\beta-1} P\left(\frac{|S_n|}{\sqrt{n \log n}} > \lambda\right) < \infty \qquad \text{for all } \lambda > 2\sqrt{\beta}$$

if

$$EY_1 = 0,\ EY_1^2 = 1 \quad \text{and} \quad E\left(\frac{|Y_1|^{4\beta+2}}{(\log(e+|Y_1|))^{2\beta+1}}\right) < \infty. \tag{3.11}$$

Combining Theorem 3.1 and Lai's [9] theorem yields the following.

COROLLARY 3.2. *Let $\beta > 0$ and suppose that*

$$(EU_{1,1})(EV_{1,1}) = 0, \qquad (EU_{1,1}^2)(EV_{1,1}^2) = 1$$

*and*

$$E\left(\frac{|U_{1,1}V_{1,2}|^{4\beta+2}}{(\log(e+|U_{1,1}V_{1,2}|))^{2\beta+1}}\right) < \infty. \tag{3.12}$$

*Then*

$$\limsup_{n\to\infty} \frac{\max_{1\le i\ne j\le n^\beta} |\sum_{k=1}^n U_{k,i}V_{k,j}|}{\sqrt{n\log n}} \le 2\sqrt{\beta} \qquad a.s.$$

PROOF. Set $a_1 = 1$, $a_n = \sqrt{n\log n}, n \ge 2$ and $p_n = [n^\beta], n \ge 1$. Then (3.1) is immediate. Note that (3.12) and (3.11) are equivalent. The first two conditions of (3.11) and Chebyshev's inequality ensure that (3.2) holds. Now (3.3) holds for all $\lambda > 2\sqrt{\beta}$ by Lai's [9] theorem. Thus by Theorem 3.1,

$$\limsup_{n\to\infty} \frac{T_n}{\sqrt{n\log n}} \le \lambda \qquad \text{a.s. for all } \lambda > 2\sqrt{\beta}.$$

The conclusion follows by letting $\lambda \downarrow 2\sqrt{\beta}$. □

Throughout the rest of this section, it is not being assumed that $\{p_n; n \ge 1\}$ is monotone.

THEOREM 3.2. *Suppose that $n/p_n$ is bounded away from 0 and $\infty$. Let $\alpha > 1/2$. Then*

$$\lim_{n\to\infty} \frac{T_n}{n^\alpha} = 0 \qquad a.s. \tag{3.13}$$

*if and only if*

$$\sum_{n=1}^\infty P\Big(\max_{1\le i\ne j\le n} |U_{1,i}V_{1,j}| \ge n^\alpha\Big) < \infty \tag{3.14}$$

*and*

$$(EU_{1,1})(EV_{1,1}) = 0 \qquad \textit{whenever } \alpha \le 1.$$

THEOREM 3.3. *Suppose that $n/p_n$ is bounded away from 0 and $\infty$. If*

$$(EU_{1,1})(EV_{1,1}) = 0, \qquad (EU_{1,1}^2)(EV_{1,1}^2) = 1$$

*and*

$$\sum_{n=1}^\infty P\Big(\max_{1\le i\ne j\le n} |U_{1,i}V_{1,j}| \ge \sqrt{n\log n}\Big) < \infty, \tag{3.15}$$

*then*

$$\limsup_{n\to\infty} \frac{T_n}{\sqrt{n\log n}} \le 2 \qquad a.s. \tag{3.16}$$

*Conversely, if*

$$\limsup_{n\to\infty} \frac{T_n}{\sqrt{n\log n}} < \infty \qquad a.s., \tag{3.17}$$

*then $(EU_{1,1})(EV_{1,1}) = 0$, $(EU_{1,1}^2)(EV_{1,1}^2) < \infty$ and (3.15) holds.*

For the proofs of Theorems 3.2 and 3.3 we need the following two lemmas.

LEMMA 3.1. *Let* $\{a_n; n \geq 1\}$ *be a nondecreasing sequence of positive constants such that*

$$\lim_{n\to\infty} a_{n+1}/a_n = 1 \quad \text{and} \quad \liminf_{n\to\infty} a_{2n}/a_n = b \in (1, \infty].$$

*Then, for every* $c > 0$ *and* $q > 1$*, the following statements are equivalent*:

$$\sum_{n=1}^{\infty} P\Big(\max_{1\leq i\neq j\leq n} |U_{1,i}V_{1,j}| \geq a_n\Big) < \infty, \tag{3.18}$$

$$\sum_{n=1}^{\infty} P\Big(\max_{1\leq i\neq j\leq cn} |U_{1,i}V_{1,j}| \geq \varepsilon a_n\Big) < \infty \qquad \textit{for all } \varepsilon > 0, \tag{3.19}$$

$$\sum_{n=1}^{\infty} P\Big(\max_{1\leq m\leq q^n} \max_{1\leq i\neq j\leq cm} |U_{m,i}V_{m,j}| \geq \varepsilon a_{[q^n]}\Big) < \infty \qquad \textit{for all } \varepsilon > 0. \tag{3.20}$$

PROOF. We only give the proof of the equivalence of (3.18) and (3.19) since the proof of the equivalence of (3.18) and (3.20) is similar. To show that (3.19) implies (3.18), note that

$$\begin{aligned} P\Big(\max_{1\leq i\neq j\leq 2cn} |U_{1,i}V_{1,j}| \geq \varepsilon a_n\Big) &\leq \sum_{l=1}^{4}\sum_{k=1}^{4} P\Big(\max_{i\neq j, i\in I_k, j\in I_l} |U_{1,i}V_{1,j}| \geq \varepsilon a_n\Big) \\ &\leq \sum_{l=1}^{4}\sum_{k=1}^{4} P\Big(\max_{1\leq i\neq j\leq cn} |U_{1,i}V_{1,j}| \geq \varepsilon a_n\Big) \\ &= 16 P\Big(\max_{1\leq i\neq j\leq cn} |U_{1,i}V_{1,j}| \geq \varepsilon a_n\Big), \end{aligned}$$

where $I_k = \{m; ((k-1)/2)cn < m \leq (k/2)cn\}$, $k = 1, 2, 3, 4$. Thus (3.19) implies that

$$\sum_{n=1}^{\infty} P\Big(\max_{1\leq i\neq j\leq 2cn} |U_{1,i}V_{1,j}| \geq \varepsilon a_n\Big) < \infty \qquad \text{for all } \varepsilon > 0. \tag{3.21}$$

Let $v$ be a positive integer such that $2^v c \geq 1$. By repeating $v-1$ times the above procedure for arriving at (3.21), we get

$$\sum_{n=1}^{\infty} P\Big(\max_{1\leq i\neq j\leq 2^v cn} |U_{1,i}V_{1,j}| \geq \varepsilon a_n\Big) < \infty \qquad \text{for all } \varepsilon > 0. \tag{3.22}$$

Thus the proof that (3.19) implies (3.18) is complete.

We now prove that (3.18) implies (3.19). Under (3.18), we can use the same idea for arriving at (3.21) to get

$$\sum_{n=1}^{\infty} P\Big(\max_{1\le i\ne j\le 2n} |U_{1,i}V_{1,j}| \ge a_n\Big) < \infty. \tag{3.23}$$

Let $b_1 \in (1, b)$. Then since $a_n \le a_{2n-1}/b_1 \le a_{2n}/b_1$ for all large $n$, it is easy to see that (3.23) implies

$$\sum_{n=1}^{\infty} P\Big(\max_{1\le i\ne j\le n} |U_{1,i}V_{1,j}| \ge (1/b_1)a_n\Big) < \infty.$$

By iterating this technique we get

$$\sum_{n=1}^{\infty} P\Big(\max_{1\le i\ne j\le n} |U_{1,i}V_{1,j}| \ge (1/b_1)^v a_n\Big) < \infty \qquad \text{for } v = 1, 2, 3, \dots. \tag{3.24}$$

Since $\lim_{v\to\infty}(1/b_1)^v = 0$, (3.19) with $c = 1$ follows from (3.24). Thus (3.22) with $c = 1$ and arbitrary $v \ge 1$ holds, and from this we get that (3.19) holds for every $c > 0$. □

LEMMA 3.2. *Suppose that $n/p_n$ is bounded away from 0 and $\infty$. Let $\{a_n; n \ge 1\}$ be as in Lemma 3.1. If*

$$\limsup_{n\to\infty} \frac{T_n}{a_n} < \infty \qquad a.s., \tag{3.25}$$

*then (3.18) holds and $(EU_{1,1})(EV_{1,1}) = 0$ whenever $\lim_{n\to\infty} a_n/n = 0$.*

PROOF. Since $n/p_n$ is bounded away from 0 and $\infty$, there exists a constant $c \ge 1$ such that $c^{-1}n \le p_n \le cn, n \ge 1$. Then it follows from (3.25) that

$$\limsup_{n\to\infty} \frac{\max_{1\le i\ne j\le c^{-1}n} |\sum_{k=1}^{n} U_{k,i}V_{k,j}|}{a_n} < \infty \qquad \text{a.s.}$$

Since $\lim_{n\to\infty} a_{n+1}/a_n = 1$,

$$\limsup_{n\to\infty} \frac{\max_{1\le i\ne j\le c^{-1}n} |\sum_{k=1}^{n+1} U_{k,i}V_{k,j}|}{a_n} < \infty \qquad \text{a.s.}$$

Note that

$$\max_{1\le i\ne j\le c^{-1}n} |U_{n+1,i}V_{n+1,j}| \le \max_{1\le i\ne j\le c^{-1}n} \left|\sum_{k=1}^{n} U_{k,i}V_{k,j}\right| + \max_{1\le i\ne j\le c^{-1}n} \left|\sum_{k=1}^{n+1} U_{k,i}V_{k,j}\right|, \qquad n \ge 1.$$

We then have

$$\limsup_{n\to\infty} \frac{\max_{1\le i\ne j\le c^{-1}n} |U_{n+1,i}V_{n+1,j}|}{a_n} < \infty \qquad \text{a.s.}$$

Then since the random vectors in the array $\{U_{k,i}, V_{k,i}; i \ge 1, k \ge 1\}$ are i.i.d., it follows from the Borel–Cantelli lemma that

$$\sum_{n=1}^{\infty} P\Big(\max_{1\le i\ne j\le c^{-1}n} |U_{1,i}V_{1,j}| \ge \lambda a_n\Big) < \infty \qquad \text{for some } \lambda > 0.$$

Using the same argument as in the proof of Lemma 3.1, we have

$$\sum_{n=1}^{\infty} P\Big(\max_{1\le i\ne j\le n} |U_{1,i}V_{1,j}| \ge \tilde{a}_n\Big) < \infty,$$

where $\tilde{a}_n = \lambda a_n$, $n \ge 1$. In view of Lemma 3.1, (3.18) follows. If $\lim_{n\to\infty} a_n/n = 0$, then (3.25) implies that

$$\lim_{n\to\infty} \frac{\sum_{k=1}^{n} U_{k,1}V_{k,2}}{n} = 0 \qquad \text{a.s.}$$

and hence by the Kolmogorov SLLN, $(EU_{1,1})(EV_{1,1}) = (EU_{1,1})(EV_{1,2}) = 0$; the proof of Lemma 3.2 is therefore complete. □

Proof of Theorem 3.2. In view of Lemma 3.2, we only need to give the proof of the "if" part. Note that (3.14) implies that

$$\sum_{n=1}^{\infty} P\Big(\max_{1\le i\le n/2} |U_{1,i}V_{1,[n/2]+i}| \ge n^{\alpha}\Big) < \infty$$

which is equivalent to

$$\sum_{n=1}^{\infty} nP(|U_{1,1}V_{1,2}| \ge n^{\alpha}) < \infty.$$

So it follows that $E|U_{1,1}V_{1,2}|^{2/\alpha} < \infty$. Setting $S_n = \sum_{k=1}^{n} U_{k,1}V_{k,2}, n \ge 1$ and applying the Baum–Katz [4] theorem, we have

$$\sum_{n=1}^{\infty} P\Big(\sup_{m\ge n} \frac{|S_m|}{m^{\alpha}} > \varepsilon\Big) < \infty \qquad \text{for all } \varepsilon > 0$$

which implies that

$$P\Big(\frac{|S_n|}{n^{\alpha}} > \varepsilon\Big) = o(n^{-1}) \qquad \text{for all } \varepsilon > 0$$

and hence by Ottaviani's inequality, it follows that

$$\max_{1\le j\le n} P\Big(\frac{|S_j|}{n^{\alpha}} > \varepsilon\Big) = o(n^{-1}) \qquad \text{for all } \varepsilon > 0. \tag{3.26}$$

Since $n/p_n$ is bounded away from 0 and $\infty$, there exists a constant $c \geq 1$ such that $c^{-1}n \leq p_n \leq cn, n \geq 1$. Thus (3.13) follows if we can show that

$$\lim_{n\to\infty} \frac{\max_{1\leq i\neq j\leq cn} |\sum_{k=1}^{n} U_{k,i}V_{k,j}|}{n^\alpha} = 0 \qquad \text{a.s.} \tag{3.27}$$

For fixed $\varepsilon > 0$, set for $1 \leq k \leq 2^n$, $1 \leq i, j \leq c2^n, n \geq 1$,

$$Y_{k,n,i,j}^{(1)} = U_{k,i}V_{k,j}I(|U_{k,i}V_{k,j}| > (\varepsilon/2^\alpha)2^{n\alpha}),$$

$$Y_{k,n,i,j}^{(2)} = U_{k,i}V_{k,j}I(|U_{k,i}V_{k,j}| \leq (\varepsilon/2^\alpha)2^{n\alpha}).$$

Then, for $2^{n-1} < m \leq 2^n, n \geq 1$,

$$\begin{aligned} \frac{\max_{1\leq i\neq j\leq cm} |\sum_{k=1}^{m} U_{k,i}V_{k,j}|}{m^\alpha} &\leq \frac{\max_{1\leq i\neq j\leq cm} |\sum_{k=1}^{m} Y_{k,n,i,j}^{(1)}|}{m^\alpha} \\ &\quad + \frac{\max_{1\leq i\neq j\leq cm} |\sum_{k=1}^{m} Y_{k,n,i,j}^{(2)}|}{m^\alpha}. \end{aligned} \tag{3.28}$$

In view of (3.14), by applying Lemma 3.1, we have

$$P\left(\max_{1\leq m\leq 2^n} \max_{1\leq i\neq j\leq cm} \left|\sum_{k=1}^{m} Y_{k,n,i,j}^{(1)}\right| = 0 \text{ eventually}\right) = 1. \tag{3.29}$$

Clearly, recalling $E|U_{1,1}V_{1,2}|^{2/\alpha} < \infty$, for all $\delta > 0$,

$$\begin{aligned} &\max_{2^{n-1}<m\leq 2^n} \max_{1\leq j\leq m} P\left(\frac{|\sum_{k=1}^{j} Y_{k,n,1,2}^{(1)}|}{m^\alpha} > \delta\right) \\ &\quad \leq 2^n P(|U_{1,1}V_{1,2}| > (\varepsilon/2^\alpha)2^{n\alpha}) \\ &\quad = o(2^{-n}) \end{aligned}$$

and this, together with (3.26), ensures that, for all $\delta > 0$,

$$\max_{2^{n-1}<m\leq 2^n} \max_{1\leq j\leq m} P\left(\frac{|\sum_{k=1}^{j} Y_{k,n,1,2}^{(2)}|}{m^\alpha} > \delta\right) = o(2^{-n}). \tag{3.30}$$

Write $\mu_{m,n} = \max_{1\leq j\leq m} \kappa_{j,m,n}$ where $\kappa_{j,m,n}$ is a median of the random variable $|\sum_{k=1}^{j} Y_{k,n,1,2}^{(2)}|/m^\alpha, 1 \leq j \leq m, 2^{n-1} < m \leq 2^n, n \geq 1$. Note that (3.30) implies that

$$\lim_{n\to\infty} \max_{2^{n-1}<m\leq 2^n} \mu_{m,n} = 0.$$

Applying Lemma 3.2 of [10] which is a generalization of the Hoffmann–Jørgensen [6] inequalities, it follows from (3.30) that for sufficiently large $n$ and every $m \in [2^{n-1}+1, 2^n]$,

$$
\begin{aligned}
&P\Big(\frac{|\sum_{k=1}^{m} Y_{k,n,1,2}^{(2)}|}{m^{\alpha}} > 10\varepsilon\Big) \\
&\quad \le P\Big(\max_{1\le k\le m} |Y_{k,n,1,2}^{(2)}| > \varepsilon m^{\alpha}\Big) + 4\Big(P\Big(\frac{|\sum_{k=1}^{m} Y_{k,n,1,2}^{(2)}|}{m^{\alpha}} > \frac{9}{2}\varepsilon - \mu_{m,n}\Big)\Big)^2 \\
&\quad \le 4\Big(P\Big(\frac{|\sum_{k=1}^{m} Y_{k,n,1,2}^{(2)}|}{m^{\alpha}} > 4\varepsilon\Big)\Big)^2 \\
&\quad \le 64\Big(P\Big(\frac{|\sum_{k=1}^{m} Y_{k,n,1,2}^{(2)}|}{m^{\alpha}} > \varepsilon\Big)\Big)^4 \\
&\quad = o(2^{-4n}).
\end{aligned}
$$

Hence

$$
\begin{aligned}
&\sum_{n=1}^{\infty} \sum_{2^{n-1}<m\le 2^n} P\Big(\frac{\max_{1\le i\ne j\le cm} |\sum_{k=1}^{m} Y_{k,n,i,j}^{(2)}|}{m^{\alpha}} > 10\varepsilon\Big) \\
&\qquad \le \sum_{n=1}^{\infty} \sum_{2^{n-1}<m\le 2^n} c^2 m^2 P\Big(\frac{|\sum_{k=1}^{m} Y_{k,n,1,2}^{(2)}|}{m^{\alpha}} > 10\varepsilon\Big) \\
&\qquad \le \sum_{n=1}^{\infty} o(2^{-n}) \\
&\qquad < \infty.
\end{aligned} \tag{3.31}
$$

Taking into account (3.29) and (3.31), we conclude from (3.28) and the Borel–Cantelli lemma that

$$
\limsup_{n\to\infty} \frac{\max_{1\le i\ne j\le cn} |\sum_{k=1}^{n} U_{k,i}V_{k,j}|}{n^{\alpha}} \le 10\varepsilon \qquad \text{a.s.}
$$

Letting $\varepsilon \downarrow 0$, (3.27) follows. □

Proof of Theorem 3.3. In view of Lemma 3.2, we only need to give the proof of the first part. Note that (3.15) implies that

$$
\sum_{n=1}^{\infty} P\Big(\max_{1\le i\le n/2} |U_{1,i}V_{1,[n/2]+i}| \ge \sqrt{n\log n}\Big) < \infty
$$

which is equivalent to

$$
\sum_{n=1}^{\infty} nP(|U_{1,1}V_{1,2}| \ge \sqrt{n\log n}\,) < \infty.
$$

So it follows that

$$E\left(\frac{U_{1,1}^4 V_{1,2}^4}{(\log(e+|U_{1,1}V_{1,2}|))^2}\right) < \infty. \tag{3.32}$$

Let $\delta > 0$ be fixed. We choose $1 < q < 2$ such that

$$\sqrt{[qn]\log[qn]} \le (1+\delta)\sqrt{n\log n} \qquad \text{for all sufficiently large } n. \tag{3.33}$$

For $1 \le k \le q^n, n \ge 1$, set

$$Z_{k,n,i,j}^{(1)} = U_{k,i}V_{k,j} I(|U_{k,i}V_{k,j}| > (\delta/2)\sqrt{q^n \log q^n}\,),$$
$$Z_{k,n,i,j}^{(2)} = U_{k,i}V_{k,j} I(q^{n/3} < |U_{k,i}V_{k,j}| \le (\delta/2)\sqrt{q^n \log q^n}\,),$$
$$Z_{k,n,i,j}^{(3)} = U_{k,i}V_{k,j} I(|U_{k,i}V_{k,j}| \le q^{n/3}).$$

Then, for $q^{n-1} < m \le q^n, n \ge 1$,

$$\begin{aligned}
&\frac{\max_{1\le i\ne j\le cm} |\sum_{k=1}^m U_{k,i}V_{k,j}|}{\sqrt{m\log m}} \\
&\quad \le \frac{\max_{1\le i\ne j\le cm} |\sum_{k=1}^m (Z_{k,n,i,j}^{(1)} - EZ_{k,n,1,2}^{(1)})|}{\sqrt{m\log m}} \\
&\qquad + \frac{\max_{1\le i\ne j\le cm} |\sum_{k=1}^m (Z_{k,n,i,j}^{(2)} - EZ_{k,n,1,2}^{(2)})|}{\sqrt{m\log m}} \\
&\qquad + \frac{\max_{1\le i\ne j\le cm} |\sum_{k=1}^m (Z_{k,n,i,j}^{(3)} - EZ_{k,n,1,2}^{(3)})|}{\sqrt{m\log m}},
\end{aligned}$$

where $c \ge 1$ is a constant such that $c^{-1}n \le p_n \le cn, n \ge 1$. Note that (3.32) ensures that

$$\max_{q^{n-1}<m\le q^n} \left|\sum_{k=1}^m EZ_{k,n,1,2}^{(1)}\right| \le \frac{4E|U_{1,1}V_{1,2}|^3}{\delta^2 \log q^n} \longrightarrow 0 \qquad \text{as } n \to \infty.$$

So, in view of Lemma 3.1, condition (3.15) implies via the Borel–Cantelli lemma that

$$\begin{aligned}
&P\left(\max_{q^{n-1}<m\le q^n} \frac{\max_{1\le i\ne j\le cm} |\sum_{k=1}^m (Z_{k,n,i,j}^{(1)} - EZ_{k,n,1,2}^{(1)})|}{\sqrt{m\log m}} \le \delta \text{ eventually}\right) \\
&\quad = 1.
\end{aligned} \tag{3.34}$$

Using the Chebyshev inequality, it follows from (3.32) that for $q^{n-1} < m \le q^n, n \ge 1$ and $\varepsilon > 0$,

$$P\left(\frac{|\sum_{k=1}^m (Z_{k,n,1,2}^{(2)} - EZ_{k,n,1,2}^{(2)})|}{\sqrt{m\log m}} > \varepsilon\right) \le \frac{mE(Z_{1,n,1,2}^{(2)})^2}{\varepsilon^2 m\log m}$$

$$\begin{aligned} &\le \frac{E|U_{1,1}V_{1,2}|^3}{\varepsilon^2 q^{n/3}\log m} \\ &= o(q^{-n/3}). \end{aligned} \tag{3.35}$$

Thus, applying Lemma 3.2 of [10] and using the same argument as in the proof of Theorem 3.2, we have that for sufficiently large $n$ and every $m \in [[q^{n-1}]+1, [q^n]]$,

$$\begin{aligned} &P\Big(\frac{|\sum_{k=1}^m (Z^{(2)}_{k,n,1,2} - EZ^{(2)}_{k,n,1,2})|}{\sqrt{m\log m}} > 46\delta\Big) \\ &\quad \le P\Big(\max_{1\le k\le m} |Z^{(2)}_{k,n,1,2} - EZ^{(2)}_{k,n,1,2}| > \delta\sqrt{m\log m}\Big) \\ &\qquad + 4\Big(P\Big(\frac{|\sum_{k=1}^m (Z^{(2)}_{k,n,1,2} - EZ^{(2)}_{k,n,1,2})|}{\sqrt{m\log m}} > 22\delta\Big)\Big)^2 \\ &\quad = 4\Big(P\Big(\frac{|\sum_{k=1}^m (Z^{(2)}_{k,n,1,2} - EZ^{(2)}_{k,n,1,2})|}{\sqrt{m\log m}} > 22\delta\Big)\Big)^2 \\ &\quad \le 4^3\Big(P\Big(\frac{|\sum_{k=1}^m (Z^{(2)}_{k,n,1,2} - EZ^{(2)}_{k,n,1,2})|}{\sqrt{m\log m}} > 10\delta\Big)\Big)^4 \\ &\quad \le 4^7\Big(P\Big(\frac{|\sum_{k=1}^m (Z^{(2)}_{k,n,1,2} - EZ^{(2)}_{k,n,1,2})|}{\sqrt{m\log m}} > 4\delta\Big)\Big)^8 \\ &\quad \le 4^{15}\Big(P\Big(\frac{|\sum_{k=1}^m (Z^{(2)}_{k,n,1,2} - EZ^{(2)}_{k,n,1,2})|}{\sqrt{m\log m}} > \delta\Big)\Big)^{16} \\ &\quad = o(q^{-(16/3)n}) \qquad [\text{by } (3.35)]. \end{aligned}$$

Hence

$$\begin{aligned} &\sum_{n=1}^\infty \sum_{q^{n-1}<m\le q^n} P\Big(\frac{\max_{1\le i\ne j\le cm} |\sum_{k=1}^m (Z^{(2)}_{k,n,i,j} - EZ^{(2)}_{k,n,1,2})|}{\sqrt{m\log m}} > 46\delta\Big) \\ &\quad \le \sum_{n=1}^\infty \sum_{q^{n-1}<m\le q^n} c^2m^2 P\Big(\frac{|\sum_{k=1}^m (Z^{(2)}_{k,n,1,2} - EZ^{(2)}_{k,n,1,2})|}{\sqrt{m\log m}} > 46\delta\Big) \\ &\quad \le \sum_{n=1}^\infty o(q^{-(7/3)n}) \\ &\quad < \infty. \end{aligned} \tag{3.36}$$

Taking into account (3.34) and (3.36), we conclude by the Borel–Cantelli lemma that

$$\limsup_{n\to\infty} \max_{q^{n-1}<m\le q^n} \frac{\max_{1\le i\ne j\le cm} |\sum_{k=1}^m H_{k,n,i,j}|}{\sqrt{m\log m}} \le 47\delta \qquad \text{a.s.,} \tag{3.37}$$

where $H_{k,n,i,j} = (Z^{(1)}_{k,n,i,j} - EZ^{(1)}_{k,n,1,2}) + (Z^{(2)}_{k,n,i,j} - EZ^{(2)}_{k,n,1,2})$. It is easy to see that $(EU^2_{1,1})(EV^2_{1,1}) = 1$ implies that for all $\varepsilon > 0$,

$$\max_{q^{n-1}<m\le q^n} P\biggl(\frac{|\sum_{k=1}^m (Z^{(3)}_{k,n,1,2} - EZ^{(3)}_{k,n,1,2})|}{\sqrt{m\log m}} > \varepsilon\biggr)$$
$$\le \frac{1}{\varepsilon^2 \log q^{n-1}} \longrightarrow 0 \qquad \text{as } n\to\infty.$$

Using the same argument as that used to obtain (3.6) in the proof of Theorem 3.1, for all large $n$

$$P\biggl(\max_{q^{n-1}<m\le q^n} \frac{\max_{1\le i\ne j\le cm} |\sum_{k=1}^m (Z^{(3)}_{k,n,i,j} - EZ^{(3)}_{k,n,1,2})|}{\sqrt{m\log m}} > 2(1+3\delta)^2\biggr)$$
$$\le 2c^2 q^{2n} P\biggl(\frac{|\sum_{k=1}^{[q^n]} (Z^{(3)}_{k,n,1,2} - EZ^{(3)}_{k,n,1,2})|}{\sqrt{[q^n]\log[q^n]}} > 2(1+2\delta)\biggr).$$

Note that

$$\sum_{k=1}^{[q^n]} E(Z^{(3)}_{k,n,1,2} - EZ^{(3)}_{k,n,1,2})^2 = [q^n] E(Z^{(3)}_{1,n,1,2} - EZ^{(3)}_{1,n,1,2})^2$$
$$\sim q^n \qquad \text{as } n\to\infty,$$
$$\max_{1\le k\le q^n} |Z^{(3)}_{k,n,1,2} - EZ^{(3)}_{k,n,1,2}| \le 2q^{n/3}, \qquad n\ge 1,$$

and

$$\lim_{n\to\infty} \frac{2(1+2\delta)\sqrt{q^n\log q^n}(2q^{n/3})}{2q^n} = 0.$$

Then, applying Lemma 7.1 of [11] which is the classical Kolmogorov exponential inequalities, we have that for all large $n$,

$$P\Biggl(\Biggl|\sum_{k=1}^{[q^n]} (Z^{(3)}_{k,n,1,2} - EZ^{(3)}_{k,n,1,2})\Biggr| > 2(1+2\delta)\sqrt{[q^n]\log[q^n]}\Biggr)$$
$$\le 2\exp\{-2(1+\delta)\log q^n\} = 2q^{-2(1+\delta)n}.$$

Hence

$$\sum_{n=1}^{\infty} P\left(\max_{q^{n-1}<m\le q^n} \frac{\max_{1\le i\ne j\le cm} |\sum_{k=1}^{m}(Z_{k,n,i,j}^{(3)} - EZ_{k,n,1,2}^{(3)})|}{\sqrt{m\log m}} > 2(1+3\delta)^2\right)$$
$$\le \sum_{n=1}^{\infty} O(q^{-2\delta n}) < \infty$$

and another application of the Borel–Cantelli lemma gives

$$\limsup_{n\to\infty} \max_{q^{n-1}<m\le q^n} \frac{\max_{1\le i\ne j\le cm} |\sum_{k=1}^{m}(Z_{k,n,i,j}^{(3)} - EZ_{k,n,1,2}^{(3)})|}{\sqrt{m\log m}} \le 2(1+3\delta)^2 \qquad \text{a.s.} \tag{3.38}$$

Combining (3.37) and (3.38) and letting $\delta\downarrow 0$, we get

$$\limsup_{n\to\infty} \frac{\max_{1\le i\ne j\le cn} |\sum_{k=1}^{n} U_{k,i}V_{k,j}|}{\sqrt{n\log n}} \le 2 \qquad \text{a.s.}$$

The proof of Theorem 3.3 is therefore complete. □

Corollary 3.3. *Let $\{X_{k,i}; k\ge 1, i\ge 1\}$ be an array of i.i.d. random variables. Suppose that $n/p_n$ is bounded away from 0 and $\infty$.*

(i) *Let $\alpha > 1/2$. Then*

$$\lim_{n\to\infty} \frac{\max_{1\le i\le p_n} |\sum_{k=1}^{n} X_{k,i}|}{n^\alpha} = 0 \qquad a.s.$$

*if and only if $E|X_{1,1}|^{2/\alpha} < \infty$ and $EX_{1,1} = 0$ whenever $\alpha\le 1$.*

(ii) *If*

$$EX_{1,1} = 0,\ EX_{1,1}^2 = 1 \quad \text{and} \quad E\left(\frac{X_{1,1}^4}{(\log(e+|X_{1,1}|))^2}\right) < \infty, \tag{3.39}$$

*then*

$$\limsup_{n\to\infty} \frac{\max_{1\le i\le p_n} |\sum_{k=1}^{n} X_{k,i}|}{\sqrt{n\log n}} \le 2 \qquad a.s. \tag{3.40}$$

*Conversely, if*

$$\limsup_{n\to\infty} \frac{\max_{1\le i\le p_n} |\sum_{k=1}^{n} X_{k,i}|}{\sqrt{n\log n}} < \infty \qquad a.s.,$$

*then*

$$EX_{1,1} = 0,\ EX_{1,1}^2 < \infty \quad \text{and} \quad E\left(\frac{X_{1,1}^4}{(\log(e+|X_{1,1}|))^2}\right) < \infty.$$

Proof. Set $U_{k,i} = X_{k,i}, V_{k,i}\equiv 1,\ k\ge 1, i\ge 1$. Then parts (i) and (ii) follow, respectively, from Theorems 3.2 and 3.3. □

**4. Proofs of the main results.** We now give the proofs of the main results.

Proof of Theorem 2.1. Theorem 2.1 follows directly from Theorem 3.2. □

Proof of Theorem 2.2. If $1/2 < \alpha \le 1$, then, in view of Remark 2.1, (2.2) implies that $EX_{1,1}^2 < \infty$. Let $\mu = EX_{1,1}$. Then $0 < \sigma^2 = E(X_{1,1} - \mu)^2 < \infty$ since $X_{1,1}$ is nondegenerate. Note that for $1 \le i \le p_n$

$$\sum_{k=1}^n (X_{k,i} - \bar{X}_i^{(n)})^2 = \sum_{k=1}^n (X_{k,i} - \mu)^2 - n(\bar{X}_i^{(n)} - \mu)^2$$

and for $1 \le i, j \le p_n$

$$\begin{aligned}\sum_{k=1}^n &(X_{k,i} - \bar{X}_i^{(n)})(X_{k,j} - \bar{X}_j^{(n)}) \\ &= \sum_{k=1}^n (X_{k,i} - \mu)(X_{k,j} - \mu) - n(\bar{X}_i^{(n)} - \mu)(\bar{X}_j^{(n)} - \mu).\end{aligned}$$

By Corollary 3.3, it is easy to see that

$$\lim_{n\to\infty} n^{1-\alpha} \max_{1\le i\le p_n} |\bar{X}_i^{(n)} - \mu| = 0 \qquad \text{a.s.},$$

$$\lim_{n\to\infty} n^{1-\alpha} \max_{1\le i<j\le p_n} |\bar{X}_i^{(n)} - \mu||\bar{X}_j^{(n)} - \mu| = 0 \qquad \text{a.s.},$$

and

$$\begin{aligned}\liminf_{n\to\infty} \min_{1\le i\le p_n} &\frac{\sum_{k=1}^n (X_{k,i} - \bar{X}_i^{(n)})^2}{n} \\ &= \liminf_{n\to\infty} \min_{1\le i\le p_n} \left( \frac{\sum_{k=1}^n (X_{k,i} - \mu)^2}{n} - (\bar{X}_i^{(n)} - \mu)^2 \right) \\ &= \liminf_{n\to\infty} \min_{1\le i\le p_n} \frac{\sum_{k=1}^n (X_{k,i} - \mu)^2}{n} \\ &\ge \liminf_{n\to\infty} \min_{1\le i\le p_n} \frac{\sum_{k=1}^n Y_{k,i}^2(b)}{n} \\ &\ge EY_{1,1}^2(b) - \limsup_{n\to\infty} \frac{\max_{1\le i\le p_n} |\sum_{k=1}^n (Y_{k,i}^2(b) - EY_{1,1}^2(b))|}{n} \\ &= EY_{1,1}^2(b) \qquad \text{a.s.},\end{aligned}$$

where $Y_{k,i}(b) = (X_{k,i} - \mu) I(|X_{k,i} - \mu| \le b), k \ge 1, i \ge 1, b > 0$. Letting $b \uparrow \infty$, we get

$$\liminf_{n\to\infty} \min_{1\le i\le p_n} \frac{\sum_{k=1}^n (X_{k,i} - \bar{X}_i^{(n)})^2}{n} \ge \sigma^2 \qquad \text{a.s.}$$

We now show that

$$\sum_{n=1}^{\infty} P\Big(\max_{1\le i<j\le n} |(X_i-\mu)(X_j-\mu)| \ge n^{\alpha}\Big) < \infty. \tag{4.1}$$

Note that

$$\begin{aligned}
&\sum_{n=1}^{\infty} P\Big(\max_{1\le i<j\le n} |(X_i-\mu)(X_j-\mu)| \ge n^{\alpha}\Big)\\
&\quad\le \sum_{n=1}^{\infty} P\Big(\max_{1\le i<j\le n} |X_iX_j| \ge \tfrac{1}{4}n^{\alpha}\Big)\\
&\qquad + 2\sum_{n=1}^{\infty} P\Big(\max_{1\le i\le n} |\mu X_i| \ge \tfrac{1}{4}n^{\alpha}\Big)\\
&\qquad + \sum_{n=1}^{\infty} P(\mu^2 \ge \tfrac{1}{4}n^{\alpha}).
\end{aligned}$$

Clearly,

$$\sum_{n=1}^{\infty} P(\mu^2 \ge \tfrac{1}{4}n^{\alpha}) < \infty.$$

By (2.2) and Lemma 3.1,

$$\sum_{n=1}^{\infty} P\Big(\max_{1\le i<j\le n} |X_iX_j| \ge \tfrac{1}{4}n^{\alpha}\Big) < \infty.$$

Also

$$\sum_{n=1}^{\infty} P\Big(\max_{1\le i\le n} |\mu X_i| \ge \tfrac{1}{4}n^{\alpha}\Big) \le \sum_{n=1}^{\infty} nP(|\mu X_1| \ge \tfrac{1}{4}n^{\alpha}) < \infty$$

since, by Remark 2.1, $E|X_1|^{2/\alpha} < \infty$. Thus (4.1) holds. Hence, applying Theorem 2.1, we have

$$\begin{aligned}
\limsup_{n\to\infty} n^{1-\alpha} L_n &\le \frac{1}{\sigma^2} \limsup_{n\to\infty} \frac{\max_{1\le i<j\le p_n} |\sum_{k=1}^{n}(X_{k,i}-\mu)(X_{k,j}-\mu)|}{n^{\alpha}}\\
&= 0 \qquad \text{a.s.}
\end{aligned}$$

This completes the proof of Theorem 2.2. □

Proof of Theorem 2.3. We first prove that (2.5) implies (2.4). Clearly, (2.4) follows from Theorem 3.3 if we can show that

$$\liminf_{n\to\infty} \frac{W_n}{\sqrt{n\log n}} \ge 2 \qquad \text{a.s.} \tag{4.2}$$

To show this, for arbitrary $b > 0$ and $k \geq 1, i \geq 1$, set

$$U_{k,i}(b) = X_{k,i} I(|X_{k,i}| \leq b) - EX_{1,1} I(|X_{1,1}| \leq b),$$
$$V_{k,i}(b) = X_{k,i} I(|X_{k,i}| > b) - EX_{1,1} I(|X_{1,1}| > b).$$

Note that

$$W_n \geq \max_{1 \leq i \neq j \leq n/c} \left| \sum_{k=1}^{n} U_{k,i}(b) U_{k,j}(b) \right| - 2 \max_{1 \leq i \neq j \leq n/c} \left| \sum_{k=1}^{n} U_{k,i}(b) V_{k,j}(b) \right| - \max_{1 \leq i \neq j \leq n/c} \left| \sum_{k=1}^{n} V_{k,i}(b) V_{k,j}(b) \right|,$$

where $c \geq 1$ is a constant such that $n/c \leq p_n \leq cn, n \geq 1$. Applying Lemma 3.1 of [7] (since $n/[n/c] \to c \in (0, \infty)$) and our Theorem 3.3, we have

$$\begin{aligned}\liminf_{n\to\infty} \frac{W_n}{\sqrt{n \log n}} &\geq \liminf_{n\to\infty} \frac{\max_{1\leq i\neq j\leq n/c} |\sum_{k=1}^{n} U_{k,i}(b) U_{k,j}(b)|}{\sqrt{n \log n}} \\ &\quad - 2 \limsup_{n\to\infty} \frac{\max_{1\leq i\neq j\leq n/c} |\sum_{k=1}^{n} U_{k,i}(b) V_{k,j}(b)|}{\sqrt{n \log n}} \\ &\quad - \limsup_{n\to\infty} \frac{\max_{1\leq i\neq j\leq n/c} |\sum_{k=1}^{n} V_{k,i}(b) V_{k,j}(b)|}{\sqrt{n \log n}} \\ &\geq 2EU_{1,1}^2(b) - 4\sqrt{EU_{1,1}^2(b)}\sqrt{EV_{1,1}^2(b)} - 2EV_{1,1}^2(b) \qquad \text{a.s.}\end{aligned}$$

Letting $b \uparrow \infty$, (4.2) follows since

$$\lim_{b\to\infty} EU_{1,1}^2(b) = 1 \quad \text{and} \quad \lim_{b\to\infty} EV_{1,1}^2(b) = 0.$$

We now show that (2.4) implies (2.5). In view of Lemma 3.2, Theorem 2.1 and Remark 2.1, (2.4) implies that $EX_1 = 0$, $EX_1^2 = \sigma^2 < \infty$ and (2.6) (with $\sqrt{n \log n}$ replaced by $\sigma^2 \sqrt{n \log n}$) holds. Hence

$$\lim_{n\to\infty} \frac{W_n}{\sqrt{n \log n}} = 2\sigma^2 \qquad \text{a.s.}$$

It follows that $2\sigma^2 = 2$ and so $EX_1^2 = 1$. Thus (2.5) holds. □

Proof of Theorem 2.4. In view of Remark 2.3, condition (2.6) implies that $EX_{1,1}^4 < \infty$. Let $\mu = EX_{1,1}$. Then $0 < \sigma^2 = E(X_{1,1} - \mu)^2 < \infty$ since $X_{1,1}$ is nondegenerate. Applying Corollary 3.3, one can see that

$$\lim_{n\to\infty} \max_{1\leq i\leq p_n} |\bar{X}_i^{(n)} - \mu| = 0 \qquad \text{a.s.},$$

$$\lim_{n\to\infty} \max_{1\leq i\leq p_n} \left| \frac{\sum_{k=1}^{n}((X_{k,i} - \mu)^2 - \sigma^2)}{n} \right| = 0 \qquad \text{a.s.}$$

and

$$\limsup_{n\to\infty}\left(\frac{n}{\log n}\right)^{1/2}\max_{1\le i\le p_n}|\bar{X}_i^{(n)}-\mu|\le 2 \qquad \text{a.s.}$$

By an argument similar to that in the proof of (4.1), (2.6) and Lemma 3.1 ensure that

$$\sum_{n=1}^{\infty}P\left(\max_{1\le i<j\le n}|(X_i-\mu)(X_j-\mu)|\ge\sigma^2\sqrt{n\log n}\right)<\infty.$$

Since, for every $i\ge 1$ and $j\ge 1$,

$$\left|\sum_{k=1}^{n}(X_{k,i}-\bar{X}_i^{(n)})(X_{k,j}-\bar{X}_j^{(n)})-\sum_{k=1}^{n}(X_{k,i}-\mu)(X_{k,j}-\mu)\right|$$
$$\le n|\bar{X}_i^{(n)}-\mu||\bar{X}_j^{(n)}-\mu|,$$

we get

$$\lim_{n\to\infty}\max_{1\le i\le p_n}\left|\frac{\sum_{k=1}^{n}(X_{k,i}-\bar{X}_i^{(n)})^2}{n}-\sigma^2\right|=0 \qquad \text{a.s.,}$$

$$\limsup_{n\to\infty}\frac{n\max_{1\le i\le p_n,1\le j\le p_n}|\bar{X}_i^{(n)}-\mu||\bar{X}_j^{(n)}-\mu|}{\sqrt{n\log n}}$$
$$=\lim_{n\to\infty}\max_{1\le i\le p_n}|\bar{X}_i^{(n)}-\mu|\times\limsup_{n\to\infty}\left(\left(\frac{n}{\log n}\right)^{1/2}\max_{1\le i\le p_n}|\bar{X}_i^{(n)}-\mu|\right)$$
$$=0 \qquad \text{a.s.,}$$

and it follows by Theorem 2.3 that

$$\lim_{n\to\infty}\left(\frac{n}{\log n}\right)^{1/2}L_n$$
$$=\frac{1}{\sigma^2}\lim_{n\to\infty}\frac{\max_{1\le i<j\le p_n}|\sum_{k=1}^{n}(X_{k,i}-\bar{X}_i^{(n)})(X_{k,j}-\bar{X}_j^{(n)})|}{\sqrt{n\log n}}$$
$$=\frac{1}{\sigma^2}\lim_{n\to\infty}\frac{\max_{1\le i<j\le p_n}|\sum_{k=1}^{n}(X_{k,i}-\mu)(X_{k,j}-\mu)|}{\sqrt{n\log n}}$$
$$=2 \qquad \text{a.s.}$$

Thus (2.11) has been established. □

**Acknowledgments.** The authors are grateful to the referee and an Associate Editor for carefully reading the manuscript and for perceptive comments which helped them improve the presentation. The authors take great pleasure to acknowledge that their work was inspired by that of Jiang [7].

Department of Mathematical Sciences
Lakehead University
Thunder Bay, Ontario
Canada P7B 5E1
e-mail: dli@lakeheadu.ca

Department of Statistics
University of Florida
Gainesville, Florida 32611
USA
e-mail: rosalsky@stat.ufl.edu